\documentclass[12pt]{article}
\usepackage{fullpage,amsfonts,amsmath,amsthm,graphicx}

\usepackage[usenames]{color}

\begin{document}

\newcommand{\mm}[1]{{\color{black}{#1}}}

\newcommand{\mmm}[1]{{\color{red}{#1}}}

\def\a{\alpha}
\def\b{\beta}
\def\c{\chi}
\def\d{\delta}
\def\D{\Delta}
\def\e{\epsilon}
\def\f{\phi}
\def\F{\Phi}
\def\g{\gamma}
\def\G{\Gamma}
\def\K{\Kappa}
\def\z{\zeta}
\def\th{\theta}
\def\Th{\Theta}
\def\la{\lambda}
\def\La{\Lambda}
\def\m{\mu}
\def\n{\nu}
\def\p{\pi}
\def\P{\Pi}
\def\r{\rho}
\def\R{\Rho}
\def\s{\sigma}
\def\S{\Sigma}
\def\t{\tau}
\def\om{\omega}
\def\Om{\Omega}
\def\smallo{{\rm o}}
\def\bigo{{\rm O}}
\def\to{\rightarrow}
\def\E{{\bf E}}
\def\ex{{\bf E}}
\def\cd{{\cal D}}
\def\rme{{\rm e}}
\def\hf{{1\over2}}
\def\R{{\bf  R}}
\def\cala{{\cal A}}
\def\cale{{\cal E}}
\def\Fscr{{\cal F}}
\def\cc{{\cal C}}
\def\calc{{\cal C}}
\def\calh{{\cal H}}
\def\call{{\cal L}}
\def\calr{{\cal R}}
\def\calb{{\cal B}}
\def\calz{{\cal Z}}
\def\calq{{\cal Q}}
\def\bk{\backslash}

\def\out{{\rm Out}}
\def\temp{{\rm Temp}}
\def\overused{{\rm Overused}}
\def\big{{\rm Big}}
\def\notbig{{\rm Notbig}}
\def\moderate{{\rm Moderate}}
\def\swappable{{\rm Swappable}}
\def\candidate{{\rm Candidate}}
\def\bad{{\rm Bad}}
\def\crit{{\rm Crit}}
\def\col{{\rm Col}}
\def\dist{{\rm dist}}

\newcommand\fix{{\rm FIX }}
\newcommand\fixx{{\rm FIX2 }}
\newcommand{\blank}{{\mathsf{Blank}}}

\newcommand{\Exp}{\mbox{\bf E}}
\newcommand{\var}{\mbox{\bf Var}}
\newcommand{\pr}{\mbox{\bf Pr}}

\newtheorem{lemma}{Lemma}
\newtheorem{theorem}[lemma]{Theorem}
\newtheorem{corollary}[lemma]{Corollary}
\newtheorem{claim}[lemma]{Claim}
\newtheorem{remark}[lemma]{Remark}
\newtheorem{observation}[lemma]{Observation}
\newtheorem{proposition}[lemma]{Proposition}
\newtheorem{definition}[lemma]{Definition}

\newcommand{\limninf}{\lim_{n \rightarrow \infty}}
\newcommand{\proofstart}{{\bf Proof\hspace{2em}}}
\newcommand{\tset}{\mbox{$\cal T$}}
\newcommand{\proofend}{\hspace*{\fill}\mbox{$\Box$}\vspace{2ex}}
\newcommand{\bfm}[1]{\mbox{\boldmath $#1$}}
\newcommand{\reals}{\mbox{\bfm{R}}}
\newcommand{\expect}{\mbox{\bf Exp}}
\newcommand{\he}{\hat{\e}}
\newcommand{\card}[1]{\mbox{$|#1|$}}
\newcommand{\rup}[1]{\mbox{$\lceil{ #1}\rceil$}}
\newcommand{\rdn}[1]{\mbox{$\lfloor{ #1}\rfloor$}}
\newcommand{\ov}[1]{\mbox{$\overline{ #1}$}}
\newcommand{\inv}[1]{\frac{1}{ #1}}

\def\calc{{\cal C}}
\def\cald{{\cal D}}

\title{The list chromatic number of graphs with small clique number}
 
\author{Michael Molloy\thanks{Dept of Computer Science,
University of Toronto, molloy@cs.toronto.edu.  Research supported by an NSERC Discovery Grant.}
 }

\maketitle

\begin{abstract} We prove that every triangle-free graph with maximum degree $\D$ has list chromatic number at most $(1+o(1))\frac{\D}{\ln \D}$. This matches the best-known upper bound for graphs of girth at least 5.  We also provide a new proof  that for any $r\geq4$ every $K_r$-free graph has list-chromatic number at most $200r\frac{\D\ln\ln\D}{\ln\D}$. 
\end{abstract}

\section{Introduction}
 We provide new proofs of two results of Johansson. The proofs are much shorter and simpler, and obtain an improvement in the constant of the first result. We use entropy compression, a powerful new take on the Lov\'asz Local Lemma. 

The first result bounds the list chromatic number of a triangle-free graph.  The list chromatic number of a graph $G$ is the smallest $q$ such that:  for any assignment of colour-lists of size $q$ to each vertex, it is possible to give each vertex  a colour from its list and obtain a proper colouring.   Johansson\cite{aj} proved that every triangle-free graph  has list-chromatic number at most $9\D/\ln\D$ where $\D$ is the maximum degree of the graph.  The leading constant was improved to 4 in~\cite{ps}.  Here we obtain $1+o(1)$:

\begin{theorem}\label{mt}
For every $\e>0$ there exists $\D_{\e}$ such that every triangle-free graph $G$ with maximum degree $\D\geq\D_{\e}$  has $\chi_{\ell}(G)\leq(1+\e)\D/\ln \D$.
\end{theorem}

In other words: every triangle-free graph with maximum degree $\D$ has list chromatic number at most $(1+o(1))\frac{\D}{\ln \D}$.

\mm{The bound in Theorem~\ref{mt}} matches the best known upper bound for graphs of girth five~\cite{jhk}, and indeed for any constant girth.  The best known lower bound is $\hf\frac{\D}{\ln \D}$ and comes from random $\D$-regular graphs.  {For constant $\D$, random $\D$-regular graphs are essentially high girth graphs: For any constant $K$, we expect $O(1)$ cycles of length greater than $K$, and so we can form a high-girth graph by removing a relatively small number of edges; furthermore, those edges form a matching, and so this changes the chromatic number by at most one.}

{This bound matches what is called the {\em shattering threshold} for colouring random regular graphs\cite{zk}, which is often  referred to as the ``algorithmic barrier''\cite{ac,zk}. 
This threshold arises in a wide class of problems on random graphs, and finding an efficient algorithm to solve any of these problems for edge-densities beyond the algorithmic barrier is a major open challenge (see e.g.~\cite{ac});  for colourings of random regular graphs, this means finding an efficient algorithm using $(1-\e)\frac{\D}{\ln \D}$ colours  for some $\e>0$.
Our proof of Theorem~\ref{mt} yields an efficient randomized algorithm to find a colouring for maximum degree up to the algorithmic barrier, not just for random regular graphs (where such algorithms are previously known~\cite{am}), but for every triangle-free graph. }

In a followup paper, Johansson\cite{aj} proved that  for any constant $r\geq 4$, every $K_r$-free graph has list-chromatic number at most $O(\D\ln\ln \D/\ln\D)$.   Here we match his bound, even when $r$ grows with $\D$.  

\begin{theorem}\label{mt2} For any  $r\geq4$, every $K_r$-free graph $G$ with maximum degree $\D$  has $\chi_{\ell}(G)\leq 200r \frac{\D\ln\ln\D}{\ln \D}$.
\end{theorem}

Theorem~\ref{mt2} holds for any $r$ but it is trivial unless $r<\ln\D/200 \ln\ln\D$.
Note also that this implies a bound on the chromatic number of $H$-free graphs for every fixed subgraph $H$, as an $H$-free graph is also $K_{|H|}$-free.  (Here $H$-free means that there is no subgraph isomorphic to $H$; the subgraph is not necessarily induced.)

These two results of Johansson were never published.  His proof for triangle-free graphs was presented in~\cite{mrbook} and his proof for $K_r$-free graphs was presented in~\cite{ps}.

 It is a longstanding conjecture\cite{aks} that for constant $r$, every $K_r$-free graph has chromatic number $O(\D/\ln\D)$.  So we make no attempt to optimize the constant in Theorem~\ref{mt2}.  Thus far, we do not even know whether the independence number is  large enough to support this conjecture.    Prior to Johansson's work, Shearer\cite{js1,js} proved that every triangle-free graph on $n$ vertices has independence number at least $(1-o(1))n\ln\D/\D$ (see also~\cite{aeks}) and that every $K_r$-free graph has independence number at least $\Omega(n\ln\D/\D\ln\ln\D)$.  His latter bound plays an important role in our proof of Theorem~\ref{mt2}.   Ajtai et al.\ conjectured that the $\ln\ln\D$ term can be removed here~\cite{aeks}.

Previous proofs of these, and similar results, used an iterative colouring procedure.  In each iteration, one would colour some subset of the vertices, where each vertex received a random colour from its list. Every vertex that received the same colour as a neighbour would be uncoloured.  (See~\cite{mrbook} for a presentation of this technique.)  One of the reasons for doing this is that the Local Lemma is much easier to apply when vertices are assigned colours independently.  Entropy compression allows us to use Local Lemma like calculations for random colouring procedures where, roughly speaking, vertices are coloured one-at-a-time with colours not appearing on any neighbours.   

This technique began with Moser's algorithm\cite{rm1} {which generated solutions to $k$-SAT whose existence was guaranteed by the Local Lemma; this was then extended by Moser and Tardos~\cite{mtl} to a very wide range of applications of the Local Lemma.} (See~\cite{tt, lf} for good expositions of the technique).  Subsequently,  Grytczuk, Kozik and Micek\cite{gkm} and Achlioptas and Iliopoulos\cite{ai} noted that this algorithm in fact can be applied to yield new existence results.  Previous applications to graph colouring (e.g.~\cite{ep, jp, ai2,pss,bcgr,djkw,gjkm}) involved situations where, throughout the algorithm, each vertex is guaranteed to have a large number of available colours to choose from.  That is not true in this paper since the degree of a vertex can be much higher than its list-size.  The novelty we use here is to treat a vertex having a small number of available colours as a bad event.

\section{Preliminary tools}
We begin with a common version of the Local Lemma; see e.g.\  Chapter 19 of \cite{mrbook}.

{\bf The Lov\'asz Local Lemma}\cite{el1} {\em Let $A_1,...,A_n$
be a set of random events, each with probability at most $\inv{4}$.  \mm{Suppose that} for each $1\leq i\leq n$ we have a subset $\cald_i$ of the  events such that $A_i$ is mutually independent
of all  other events outside of $\cald_i$. \mm{If} for each $1\leq i\leq n$ we have
\[\sum_{j\in\cald_i}\pr(A_j)<\inv{4},\]
then $\pr(\ov{A_1}{\cap...\cap}\ov{A_n})>0$.
}
\vskip0.2cm

We say that boolean variables $X_1,...,X_m$ are {\em negatively correlated} if
\[\mbox{ for all } I\subseteq\{1,...,m\}: \hspace{5ex} \pr\left(\wedge_{i\in I}X_i\right)\leq \prod_{i\in I}\pr(X_i).\]
Panconesi and  Srinivasan\cite{ps2} noted that many Chernoff-type bounds on independent variables also hold on negatively correlated variables.  We will use the following:

{
\begin{lemma}\label{lncv}
Suppose $X_1,...,X_m$ are  boolean variables, and set $Y_i=1-X_i$.
Set $X=\sum_{i=1}^m X_i$.  Then for any $0<t\leq \Exp(X)$:
\begin{enumerate}
\item[(a)] If $X_1,...,X_m$ are negatively correlated then
$\pr(X>\Exp(X)+t)<e^{-t^2/3E(X)}$.
\item[(b)] If $Y_1,...,Y_m$ are negatively correlated then
$\pr(X<\Exp(X)-t)<e^{-t^2/2E(X)}$.
\end{enumerate}
\end{lemma}

In this paper, we only require part (b).

Part (a) follows from Corollary 3.3 of~\cite{ps2}. The proof of part (b) is very similar and we sketch it here.

For independent variables, the bound follows from {standard Chernoff-type bounds; e.g.\ we refer to} Theorem 2.3(c) in~\cite{cmsurv}.  To adapt the proof so that it holds when $Y_1,...,Y_m$ are negatively correlated, we only need one change. Set $Y=\sum_{i=1}^m Y_i = m-X$. The proof for independent variables uses that for  any $h>0$:
\[\Exp(e^{hY})=\Exp(\prod_{i=1}^me^{hY_i})= \prod_{i=1}^m\Exp(e^{hY_i}).\]
We replace this with 
\begin{equation}\label{ehy}
\Exp(e^{hY})=\Exp(\prod_{i=1}^me^{hY_i})\leq \prod_{i=1}^m\Exp(e^{hY_i}).
\end{equation}
The highlights of the proof from~\cite{cmsurv} are:  Set $p_i=\pr(Y_i)$ for each $i$ and set $p=\sum p_i/m=\ex(Y)/m$.  For any $h>0$ we have $\Exp(e^{hY_i})=1-p_i+p_ie^h$ and so~(\ref{ehy}) and the arithmetric mean-geometric mean inequality yield
\[\Exp(e^{hY})\leq \prod_{i=1}^m(1-p_i+p_ie^h)\leq(1-p+pe^h)^m.\]
Thus $\pr(Y\geq s)\leq e^{-hs}(1-p+pe^h)^m$.  A good choice of $h$ (see the proof of Lemma 2.2 in~\cite{cmsurv}) yields that for any $0\leq z\leq 1$,
\[\pr(X\leq E(X)-mz)=\pr(Y\geq E(Y)+mz)\leq\left(\left(\frac{p}{p+z}\right)^{p+z}\left(\frac{1-p}{1-p-z}\right)^{1-p-z}\right)^m.\]  Now set $t=mz$ and apply some calculus (see the proof of  Lemma 2.3(c) in~\cite{cmsurv}) to obtain the bound for Lemma~\ref{lncv}(b).

{\bf Remark:}  Intuitively, it seems that when $X_1,....,X_m$ are negatively correlated then typically $Y_1,...,Y_m$ would also be negatively correlated.  Indeed that is the case in the application of Lemma~\ref{lncv} in this paper.   However, it is not always the case.   Choose a string from an urn containing two copies of the strings $\{000,011,101,110\}$ and one copy of each of the other boolean strings of length three.  Let $X_i$ be the event that the $i$th digit is 1.  Then $X_1,X_2,X_3$ are negatively correlated but\mm{ $\pr(Y_1\wedge Y_2\wedge Y_3)=\frac{1}{6}>\frac{1}{8}=\pr(Y_1)\pr(Y_2)\pr(Y_3)$.
}

\section{Triangle-free graphs}\label{stf}

{Each vertex $v$ has a list $\calc_v$ of colours  that may be assigned to $v$} of size 
\[|\calc_v|=q:=(1+\e)\frac{\D}{\ln\D}.\]
It suffices to prove Theorem~\ref{mt} for small $\e$;  in particular we will assume $\e<1$.

 A partial list colouring $\s$ is a colour assignment to a subset of the vertices, where the colours are drawn from their lists.  Given a partial colouring, {it is helpful if each vertex has many colours which do not appear on its neighbourhood.  To this end,} we set 
\[L=\D^{\e/2}.\]

Note that if $\D$ neighbours of $v$ are each independently given a uniformly random colour from their lists, then the expected number of colours from $\calc_v$ that are not chosen for any neighbour of $v$ is at least $q\left(1-1/q \right)^{\D}\approx (1+\e)\D^{\frac{\e}{1+\e}}/\ln\D>L$. So it is plausible that we can obtain a colouring in which every vertex has {at least $L$ colours which do not appear on its neighbourhood}. In fact we will prove that we can obtain {such a {\em partial} colouring with a substantial number of vertices coloured.}  From this, it will be straightforward to complete the colouring.

It will be convenient to treat $\blank$ as a colour, and  the uncoloured vertices are viewed as having been assigned this colour. $\blank$ is the only colour that can be assigned to two neighbours.
Most of our work goes towards finding a partial list colouring with certain properties that make it easy to complete to a full colouring. 

We use $N_v$ to denote the open neighbourhood of $v$ (to be clear:  $v\notin N_v$.)
Given a partial colouring $\s$, we define for each vertex $v$ and colour $c\neq\blank$:

\begin{eqnarray*}
&&L_v \mbox{ is the set of colours in $\calc_v$ not appearing on $N_v$, along with } \blank;
\\
&&T_{v,c} \mbox{ is the set of vertices $u\in N_v$ such that $\s(u)=\blank$ and $c\in L_u$}.\\
\end{eqnarray*}
Note that the preceding definition does not apply to $T_{v,\blank}$; it will be convenient to set $T_{v,\blank}=\emptyset$ for all $v$.
  
Given a partial colouring, we define the following two {\em flaws} for any vertex $v$:

\begin{eqnarray*}
B_v&\equiv & |L_v|<L\\
Z_v &\equiv& \sum_{c\in L_v}|T_{v,c}| >\inv{10} L\times |L_v|
\end{eqnarray*}

We say $v$ is the {\em vertex} of flaw $f=B_v$ or $Z_v$, and we denote $v(f):=v$.

{
\begin{observation}\label{odist} $B_v$ is determined by the colours of the vertices in $N(v)$ and $Z_v$ is determined by the colours of the vertices within distance two of $v$.
\end{observation}
}

{\bf Remark} If we were content with proving the weaker bound of $\chi_{\ell}(G)<(2+o(1))\frac{\D}{\ln\D}$ colours, then we could have defined $Z_v$ to be a much simpler flaw, namely that $v$ has at least $L$ blank neighbours.  We use that flaw in Section~\ref{skr}.

Our main goal is to find a partial colouring which has no flaws.  The following proof that such a colouring can be completed to a proper colouring with no blank vertices is essentially the proof of the main result in~\cite{brlist}.

\begin{lemma}\label{lfinish}
Suppose we have a  partial list colouring $\s$ such that for every vertex $v$, neither $B_v$ nor $Z_v$ hold.  Then we can colour the blank vertices to obtain a full list colouring.
\end{lemma}

\proofstart  We give each blank vertex $v$ a uniformly chosen colour from $L_v\bk\blank$.  For any edge $uv$ and colour  $c\in L_u\cap L_v$, $c\neq\blank$ we define $A_{uv,c}$ to be the event that $u,v$ both receive $c$.  Then $\pr(A_{uv,c})= 1/(|L_u|-1)(|L_v|-1)$.  Furthermore, $A_{uv,c}$ shares a vertex with at most $\sum_{c'\in L_v}|T_{v,c'}|+\sum_{c'\in L_u}|T_{u,c'}|$ other events.  The number of such events is at most $\inv{10}L(|L_v|+|L_u|)$ since $Z_u,Z_v$ do not hold. It is straightforward to check that $A_{uv,c}$ is mutually independent of all events with which it does not share a vertex (see e.g.\ the Mutual Independence Principle in Chapter 4 of~\cite{mrbook}).
So our lemma follows from the Local Lemma as $B_u,B_v$ do not hold and so
\begin{eqnarray*}
\inv{(|L_u|-1)(|L_v|-1)}\times\frac{L(|L_v|+|L_u|)}{10}&\leq&{\frac{L}{10(|L_u|-1)}\times\frac{|L_v|}{|L_v|-1}
+ \frac{L}{10(|L_v|-1)}\times\frac{|L_u|}{|L_u|-1}}\\
&<&{\inv{9}+\inv{9}}<\inv{4},
\end{eqnarray*}
for $\D>20^{2/\e}$; i.e.\ $L>20$.
\proofend

In the next section, we will present an algorithm to find a flaw-free colouring.

\subsection{Our colouring algorithm}\label{sca}

Consider a partial colouring $\s$ and any flaw $f$ of $\s$.  We will use a recursive algorithm to correct $f$.
Recall that every neighbourhood is an independent set, and so we recolour the vertices in a neighbourhood independently. 

We use the following ordering on the flaws: every $B_v$ comes before every $Z_u$, and the $B_v$'s and $Z_u$'s are each ordered according to the labels of $v,u$.
We use $\dist(w,v)$ to denote the distance from $w$ to $v$; i.e. the number of edges in a shortest $w,v$-path.
\begin{tabbing}
{\bf FIX}($f,\s$)\\
Set  \= $v=v(f)$ and assign each $u\in N_v$ a uniformly selected colour from $L_u$. \\
While there are any flaws $B_w$ with $\dist(w,v)\leq 2$ or $Z_w$ with $\dist(w,v)\leq 3$:\\
\>Let $g$ be the least such flaw and call \fix($g,\s'$) where $\s'$ is the current colouring.\\
Return the current colouring.
\end{tabbing}

\noindent{\bf Remark} It is possible that $f$ still holds after recolouring the neighbourhood of $f$, but then $f$ itself would count as a flaw within distance  2 or 3 in the next line (but is not necessarily the least of those flaws).  Note further that even if $f$ does not hold after the recolouring, it is possible for future recolourings to bring $f$ back and so {FIX} may be called again on $f$ further down in the recursive calls.  


%


Next we note that if FIX terminates, then we have made progress in correcting  the flaws.

\begin{observation}\label{ofix} In the colouring returned by {FIX}($f,\s$):
\begin{enumerate} 
\item[(a)]  $f$ does not hold; and
\item[(b)] there are no flaws that did not hold in $\s$.
\end{enumerate}
\end{observation}

\proofstart Part (a) is true because we cannot exit the while loop if $f$ holds.  Part (b) is true because any new flaw $f'$ must have arisen  {during a call of  \fix on some $f''$ whose vertex is within distance two or three of $v(f')$ (depending on whether $f'$ is a $B$-flaw or a $Z$-flaw), as these are the only calls in which  a vertex within distance one or two of $v(f')$ can be recoloured (see Observation~\ref{odist}).} But we would not have exited the while loop of that call if $f'$ still held.
\proofend

So we can obtain a flaw-free colouring by beginning with any partial colouring, e.g.\ the all-blank colouring, and then calling \fix at most once for each of the at most $2n$ flaws of that colouring.  Thus it suffices to prove that \fix terminates with positive probability; in fact, we will show that with high probability it terminates quickly (see the remark at the end of Subsection~\ref{sat}).

In the next subsection we prove that the proportion of colourings of  $N(v)$ for which $f$ holds is at most $\D^{-4}$.  
In Subsection~\ref{sat} we use that to show  \fix terminates.   Note that there are at most $2\D^3$ flaws $g$ which could appear in the the while loop in \fix$(f,\s)$.  Since $2\D^3\times \D^{-4}<\inv{4}$ (for large $\D$) this feels like a Local Lemma computation.  Entropy compression allows us to use such a computation in a  procedure like FIX, which is more complicated than what we would typically apply the Local Lemma to; in particular note how quickly dependency spreads amongst the various flaws while running FIX.

\subsection{Probability bounds}\label{spb}
In this section, we prove the key bounds on the probability of our flaws.  

\noindent{\bf Setup for Lemma~\ref{lbv}:}  Each vertex $u\in N_v$ has a list $L_u$ containing $\blank$ and perhaps other colours.  We give each $u\in N_v$ a random colour from $L_u$, where the choices are made independently and uniformly. This assignment determines $L_v, T_{v,c}$.

\begin{lemma}\label{lbv}    
\begin{enumerate}
\item[(a)] $\pr(|L_v|<L)<\D^{-4}$.
\item[(b)] $\pr( \sum_{c\in L_v}|T_{v,c}| >\inv{10} L\times |L_v|)<\D^{-4}$.
\end{enumerate}
\end{lemma}

{\bf Remarks}\\
\noindent (1) This looks like an analysis of the probability that the recolouring in \mm{the first line of \fix} produces another flaw on $N_v$.  But we will actually apply it to count the number of choices for the flawed colouring that was on $N_v$ \mm{before the recolouring}. This subtlety is important if one attempts to adapt this proof by using a different recolouring procedure designed to have a low probability of producing a flaw. \\
\noindent (2) Kim's proof~\cite{jhk} for graphs of girth five was much simpler than Johansson's proof~\cite{aj} for triangle-free graphs.  The main reason was that if $G$ has girth five then the neighbours of $v$ have disjoint neighbourhoods (other than $v$) which resulted in their lists being, in some sense, independent of each other.  In a triangle-free graph with many 4-cycles, we could have two neighbours $u_1,u_2$ of $v$ whose neighbourhoods overlap a great deal and thus their lists would be highly dependent.  Intuitively, it was clear that this should be helpful: if $L_{u_1}$ and $L_{u_2}$ are very similar then $u_1,u_2$ would tend to get the same colour which would tend to increase the size of $L_v$.  But, frustratingly, we did not know how to take advantage of this.   In the current paper, the fact that dependencies between $L_{u_1},L_{u_2}$ do not hurt is captured by the {stronger fact that Lemma~\ref{lbv} holds for any set of lists on the neighbours of $v$, even lists produced by an adversary.}

\proofstart
For each colour $c\in\calc_v\bk \{\blank\}$ we define:
\[\r(c)=\sum_{u\in N_v: c\in L_u}\inv{|L_u|-1}.\]
Thus, since each $L_u$ has $|L_u|-1$ non-Blank colours,
\begin{equation}\label{eem}
\sum_{c\in\calc_v\bk \{\blank\}} \r(c)\leq \sum_{u\in N_v}\sum_{c\in L_u\bk \{\blank\}}\inv{|L_u|-1}\leq \D.
\end{equation}
{\em Part (a):}  If  $c\in L_u$ then $|L_u|\geq 2$ and so we have $1-\inv{|L_u|}>e^{-1/(|L_u|-1)}$. We apply this inequality to obtain:
\begin{equation}\label{eel}
\ex(|L_v|)=1+\sum_{c\in\calc_v\bk  \{\blank\}}\prod_{u\in N_v: c\in L_u}
\left(1-\inv{|L_u|}\right)>\sum_{c\in\calc_v\bk  \{\blank\}}e^{-\r(c)}.
\end{equation}
By convexity of $e^{-x}$, (\ref{eem}) and recalling that $|\calc_v|=q=(1+\e)\D/\ln \D$ we have
\[\ex(|L_v|)>q e^{-\D/q}=\frac{(1+\e)\D}{\ln\D}\times\D^{-\frac{1}{1+\e}}>2\D^{\e/2}=2L,\]
for $\e<1$.

To prove concentration, we set $X_c$ to be the indicator variable that $c\in L_v$; thus  $|L_v|=1+\sum_{c\in\calc_v\bk \{\blank\}}X_c$.   {We wish to apply Lemma~\ref{lncv}(b) to bound the probability that $|L_v|$ is too small, and so we set $Y_c=1-X_c$ and argue that the variables $\{Y_c\}$ are negatively correlated.

{\em Claim:  For any $I\subseteq \calc_v\bk \{\blank\}$, $\pr(\wedge_{c\in I}Y_c)\leq\prod_{c\in I}\pr(Y_c)$.}

{\em Proof:} Consider  any $I\subseteq \calc_v\bk \{\blank\}$ and $c'\notin I$. We will first argue that 
\begin{equation}\label{eyc1}
\pr( \wedge_{c\in I}Y_c| X_{c'} )\geq \pr(\wedge_{c\in I}Y_c).
\end{equation}
To sample a colour assignment conditional on $X_{c'}$ we simply choose for each  $u\in N_v$, a uniform colour from $L_u\bk \{c'\}$.  Since $c'\notin I$, it is clear that this does not decrease the probability that every colour in $I$ is selected at least once, i.e. $\pr(\wedge_{c\in I}Y_c)$.   This establishes~(\ref{eyc1}).  This is equivalent to $\pr( \wedge_{c\in I}Y_c| Y_{c'} )\leq \pr(\wedge_{c\in I}Y_c)$, which is equivalent to
\begin{equation}\label{eyc2}
\pr( Y_{c'}|\wedge_{c\in I}Y_c ) \leq \pr(Y_{c'}).
\end{equation}
Applying~(\ref{eyc2}) inductively yields the claim.
\proofend
}

 Now Lemma~\ref{lncv}(b) yields:
\[\pr(|L_v|<\hf\ex(|L_v|))< e^{-\inv{8}\ex(|L_v|)}<e^{-\inv{4}\D^{\e/2}}<\D^{-4},\]
for $\D$ sufficiently large in terms of $\e$.  This proves part (a).

{\em Part (b):}  Let $\Psi$ be the set of colours $c\in L_v\bk \{\blank\}$ with $\r(c)>\D^{\e/4}$.  Using the same calculations as those for~(\ref{eel}), but this time applying $1-\inv{|L_u|}<e^{-1/|L_u|}<e^{-1/2(|L_u|-1)}$ for $|L_u|\geq 2$, the probability that $L_v$ contains at least one colour from $\Psi$ is at most
\[\ex(|L_v\cap\Psi|)<\sum_{c\in\Psi} e^{-\hf\r(c)}<q e^{-\hf\D^{\e/4}}<\hf\D^{-4},\]
for $\D$ sufficiently large  in terms of $\e$.  For any $c\notin\Psi$:
\[\ex(|T_{v,c}|)=\sum_{u:c\in L_u}\inv{|L_u|}<\r(c)\leq\D^{\e/4}.\]
Since the choices of whether $u\in T_{v,c}$, i.e. whether $u$ receives $\blank$, are made independently, {standard concentration bounds apply.  E.g.\ Theorem 2.3(b) of~\cite{cmsurv} says that for any $\e>0$,
\[\pr(|T_{v,c}|>(1+\e)\ex(|T_{v,c}|)<e^{-\e^2\ex(|T_{v,c}|)/2(1+\frac{\e}{3})},\] 
which yields } $\pr(|T_{v,c}|>\ex(|T_{v,c}|)+\D^{\e/4})<e^{-\frac{3}{8}\D^{\e/4}}$. So
the probability that there is at least one $c\notin\Psi$ with $|T_{v,c}|>2\D^{\e/4}$ is at most
\[qe^{-\frac{3}{8}\D^{\e/4}}<\hf\D^{-4},\]
for sufficiently large $\D$.  So with probability at least $1-\D^{-5}$ we have
\[\sum_{c\in L_v\bk  \{\blank\}}|T_{v,c}| = \sum_{c\in L_v\bk\Psi}|T_{v,c}|\leq 2\D^{\e/4}|L_v|<\inv{10}L\times |L_v|.\]
\proofend


\subsection{The algorithm terminates}\label{sat}
 The basic idea behind entropy compression is that a string of random bits cannot be represented by a shorter string.  We will consider the string of random bits used for the recolouring steps of \fix and show that as we run \fix we can record a file which allows us to recover those random bits.  Each time we call \fix($g,\s$), we record the name of $g$ and the colours of the vertices that determine $g$.  It is not hard to see that this, along with the current colouring, will allow us to reconstruct all of the preceding  random colour choices.  Because the colours which determine $g$ indicate that something unlikely occurred (namely the flaw $g$), we can represent those colours in a very concise way.  However, it may take a large amount of space to record the name of $g$. So instead, we use the degree bound in our graph to record a concise piece of information that will allow us to determine the name of $g$. This will lead to a compression of those random colour choices if the algorithm continues for too many steps.

First we describe these concise representations.  Consider any vertex $v$.  Let $N^3(v)$ denote the set of vertices within distance 3 of $v$ (including $v$ itself).  For each $1\leq \ell\leq |N^3(v)|<\D^3$ we let $\omega(\ell,v)$ denote the $\ell$th vertex of $N^3(v)$ when those vertices are listed in order of their labels.   When we call, e.g. \fix($B_w,\s'$) while running \fix($Z_v,\s''$), rather than recording the name ``$B_w$'' it will suffice to just record ``$(B,\ell)$'' where $w=\omega(\ell,v)$.  So despite the fact that the number of vertices, and hence the size of the label of $w$, is not bounded in terms of $\D$, we are able to record $w$ using only roughly $3\log_2\D$ bits.

Suppose that we are given a collection of lists $\call=\{L_u:u\in N_v\}$ of available colours for the neighbours of $v$. Let $\calb(\call)$, resp. $\calz(\call)$ be the set of all colour assignments from these lists such that $B_v$, resp. $Z_v$, holds.
Lemma~\ref{lbv} implies that $|\calb(\call)|,|\calz(\call)|<\D^{-4}\prod_{u\in N_v}|L_u|$.  For each $1\leq\ell\leq |\calb(\call)| +|\calz(\call)|$, we let $\beta(\ell,\call)$ denote the $\ell$th member of $\calb(\call)\cup\calz(\call)$ in some fixed ordering. When we run, e.g.\  \fix($B_v$) we record the colours of $N_v$ {\em before they get recoloured}; but instead of listing all the colours, we only need to record the value $\ell$ such that those colours are $\b(\ell,\call)$.

We add some write statements to \fix as follows.

\begin{tabbing}
{\bf FIX}($f,\s$)\\
Set \= $\call=\{L_u:u\in N_{v(f)}\}$.\\
\>  {\em Write} ``COLOURS =  $\ell$'' where {$\beta(\ell,\call)$ is the colouring of $N_{v(f)}$}.\\
(*) Set  \= $v=v(f)$ and assign each $u\in N_v$ a uniformly selected colour from $L_u$.\\
While there are any flaws $B_w$ with $\dist(w,v)\leq 2$ or $Z_w$ with $\dist(w,v)\leq 3$:\\
\>Let $g$ be \=the least such flaw and call {\bf FIX}($g,\s'$) where $\s'$ is the current colouring.\\
\>\>{\em Write}  ``FIX (B,$\ell$)'' or \=``FIX (Z,$\ell$)'' (depending on whether $g$ is a B-flaw or an Z-flaw)\\
\>\> \>where $v(g)=\omega(v(f),\ell)$\\
Return the current colouring.\\
\>\>{\em Write }``Return''
\end{tabbing}

Let $\s_0$ be any initial colouring and let $f$ be any flaw of $\s_0$.  We will analyze a run of \fix($\s_0,f$).
After $t$ executions of the line (*) we set
\begin{eqnarray*}
&& \s_t  \mbox{ is the current colouring}\\
&& H_t  \mbox{ is the file that we write to}\\
&& R_t  \mbox{ is the string of random bits that were used for all executions of (*)}
\end{eqnarray*}

In our formal proofs, we will not in fact make use of $R_t$; we only use it to give an intuitive picture of the compression of our random bits.  Thus we are not careful about issues such as ensuring that each random choice uses an integer number of bits.

\begin{lemma}\label{lrecon}  Given $\s_0,\s_t, f, H_t$ we can reconstruct the first $t$ steps of FIX.
\end{lemma}

\proofstart Let $f_i$ denote the flaw addressed during the $i$th execution of (*).   First observe that $f_1,...,f_t$ can be determined by $\s_0, f, H_t$.  Indeed,  proceed inductively: We know the sequence $f_1=f,...,f_{i-1}$.  \fix($f_i,\s_{i-1}$) was called while executing \fix($f_j,\s_{j-1}$) for some $j<i$. The locations of the ``Return'' lines in $H_t$ are enough to determine the value of $j$, and by induction we know $f_j$. So the $i$th ``FIX (-,$\ell$)'' line tells us that $v(f_i) = \omega(v(f_j),\ell)$ and also tells us whether $f_i=B_{v(f_i)}$ or $Z_{v(f_i)}$.  

Next observe that, having determined $f_1,...,f_t$, we can reconstruct the colours assigned in each execution of (*) from $H_t$ and $\s_t$.  To see this, note that we can reconstruct $\s_{t-1}$ from $H_t,\s_t,f_t$. We know that $\s_{t-1}=\s_{t}$ on all vertices other than $N_{v(f_t)}$.  This and the fact that our graph is triangle-free imply that  for every $u\in N_{v(f_t)}$, the list $L_u$ does not change during step $t$.   {So the collection of lists  $\call=\{L_u:u\in N(v(f_t))\}$ does not change during the $t$th recolouring and so $\s_t$ and the $t$th ``COLOURS=$\ell$'' line allows us to recover $\s_{t-1}(N_{v(f_t)})=\beta(\ell,\call)$.}
Furthermore, $\call$ and $\s_t(N_v)$ tell us what colours were selected during the $t$th execution of (*).  Working backwards, this determines $\s_t,\s_{t-1},...,\s_1$ and hence all of our random choices.
\proofend

So $R_t$ can be represented by $(\s_0,\s_t, f, H_t)$.  {The essence of the remainder of our argument is } that if \fix($\s_0,f$) continues for $t$ steps, where $t$ is large, then $(\s_0,\s_t,f,H_t)$ when expressed in binary will be much shorter than $R_t$.  Any method to represent a random string of bits by a much shorter string must fail w.h.p.   So this implies that w.h.p.\ we  terminate before very many steps.

{The rough idea is:   During the $i$th execution of (*), recall that $f_i$ is the  flaw being addressed and define:
\[ \La_i =\prod_{u\in N_{v(f_i)}} |L_u| \mbox{ at the time of the $i$th execution of (*). }  \]
The $i$th execution of (*) selects one of $\La_i$ possible colourings of $N_{v(f_i)}$ and so the total number of random bits used during the first $t$ executions is $\sum_{i=1}^t \log_2\La_i$.  Note that this number depends on the actual random choices that are made.

After $t$ executions of (*)  $H_t$ consists of: (a)  $t-1$ ``FIX(-,$\ell$)'' lines in which $\ell<\D^3$; (b)  $t$ ``COLOURS =  $\ell$'' lines in which the $i$th such line has $\ell\leq|\calb(\call)| +|\calz(\call)|  <2\D^{-4}\La_i$; (c)  fewer than $t$ ``Return'' lines.  So the total number of bits required to record $H_t$ is 
\[\sum_{i=1}^t [3\log_2\D+\log_2(2\D^{-4}\La_i)+O(1)]=-t(\log_2\D+O(1))+\sum_{i=1}^t \log_2\La_i.\]  
Thus in each execution of (*) writing to $H_t$ requires roughly $\log_2\D$ fewer bits than the number of random bits added to $R_t$.

Letting $n$ be the number of vertices, the number of  choices for each of the partial list colourings $\s_0,\s_t$ is at most $q^n$ and there are $2n$ choices for $f$.  So to record $(\s_0,\s_t,f)$ requires $2n\log_2 q + \log_2 n +1<2n\log_2\D$ bits (for sufficiently large $n$). The main point is that this does not change with $t$ and so if $t$ is large in terms of $n$ then  $|(\s_0,\s_t,f,H_t)|\leq |R_t| $, as required. 

Annoying technical issues arise when $\La_i$ is not a power of 2, and so our formal proof will use direct  probability bounds in which sizes of  the bitstreams are only implicit. 
} 

\begin{lemma}\label{cshort}
For any partial colouring $\s$ and any flaw $f$ of $\s$, the probability that \fix$(f,\s)$ continues for at least $2n$ executions of (*) is at most $\D^{-n/2}$, \mm{where $n$ is the number of vertices}.
\end{lemma}

\proofstart Set $T=2n$ and run  \fix$(f,\s)$ until it either terminates or carries out $T$ executions of (*).  

Let $\calq$ be any possible run of  \fix$(f,\s)$ that lasts for at least $T$ executions.   At the $i$th execution, recall that $\La_i=\prod_{u\in N_{v(f_i)}} |L_u|$ is the number of choices for the recolouring.  We choose this  recolouring by taking a uniform integer $x_i$ from $\{1,...,\La_i\}$.  Note that $\La_i$ is determined by $f,\s$ and $x_1,...,x_{i-1}$.  Set $\La=\La(\calq)=\prod_{i=1}^T{\La_i}$ and set $\la=\la(\calq)=\rdn{\log_2 \La}$ {(intuitively, $\la$ can be thought of as the number of random bits generated).} The probability that we carry out the run $\calq$ is $1/\La\leq 2^{-\la}$. 

 Note that $\La_i\leq(q+1)^\D$ for each $i$ and so $\la<T\D\log_2(q+1)<T\D\log\D$.

{
Given $\s_0=\s$ and $f$, Lemma~\ref{lrecon} says that $H_T,\s_T$ determine $\calq$.  So we will enumerate the number of choices for $\calq$ by enumerating the number of choices for $(H_T,\s_T)$.  We will do this by considering the size of a string encoding $(H_T,\s_T)$ in binary.

The number of choices for $\s_T$ is $(q+1)^n$, so it can be recorded with $\rup{n\log_2(q+1)}$ bits.  The $i$th line of $H_T$ consists of:  (1) a FIX line containing a number of size at most $2\D^3$; it requires $3\log_2 \D+O(1)$ bits; (2) we either do or do not write a ``Return'' line; this costs $O(1)$ bits; (3) a COLOURS line containing a number of size at most  $2\D^{-4}\La_i$; it requires $\log_2 \La_i-4\log_2\D+O(1)$ bits.  So the total size of the string recording $(H_T,\s_T)$ and hence recording $\calq$ is at most 
\[n\log_2(q+1)+\log_2\La(\calq)-T(\log_2\D - O(1))<\la(\calq)-\frac{2}{3} n\log_2 \D,\]
for $\D$ sufficiently large and \mm{since} $T=2n$.   So the total number of choices for a run $\calq$ of length $T$ and with $\la(\calq)=\la$ is at most $2^{\la-\frac{2}{3} n\log_2 \D}=2^{\la}\D^{-2n/3}$.  Thus the probability that we continue for $T=2n$ steps is at most
\[\sum_{\la= 1}^{T\D\log\D} 2^{-\la}\times 2^{\la}\D^{-2n/3}=2n\D\log\D\times\D^{-2n/3}<\D^{-n/2}.\]
\proofend
}

\subsection{Proof of Theorem~\ref{mt}}

As described above, the results of the preceding subsections provide a proof of Theorem~\ref{mt}:

{\bf Proof of Theorem~\ref{mt}} Consider any $\e>0$ and any assignment of lists of size $q=(1+\e)\D/\ln\D$ colours to the vertices.
We begin by assigning $\blank$ to every vertex. Then we repeatedly call \fix to eliminate any remaining flaws.  More formally:  While there is any flaw $f$ we call \fix($f,\s$) where $\s$ is the current partial colouring.  By Lemma~\ref{cshort} each call terminates within $O(n)$ executions of (*) with probability at least $1-\D^{-n/2}$. By Observation~\ref{ofix}, the number of flaws decreases by at least one after each call. There are at most $2n$ initial flaws and so we obtain a flaw-free partial colouring $\s^*$ after at most $2n$ calls of \fix($f,\s$) with probability at least $1-2n\D^{-n/2}>0$.  Lemma~\ref{lfinish} implies that the $\blank$ vertices of $\s^*$ can be recoloured to give the required proper list colouring.
\proofend

{\bf Remark} This easily yields a polytime algorithm to produce the list colouring. Calling \fix at most $2n$ times w.h.p.\ produces $\s^*$ in $O(n^2\D^2\ln\D)$ time; in fact, extending the definition of $H_t,R_t,\s_t$ to cover the sequence of colourings/executions produced over the sequence of at most $2n$ calls of \fix can reduce this running time to $O(n\ln n\D^2\ln\D)$ (see e.g.\ the approach in~\cite{ai}).  The main result of~\cite{mtl} yields a polytime algorithm corresponding to Lemma~\ref{lfinish}, which we use to complete the colouring.

\section{$K_r$-free graphs}\label{skr}
With a more complicated recolouring step, the same proof can be adapted to $K_r$-free graphs.  The setup is the same as in Section~\ref{stf} except with a larger list size:

{Each vertex $v$ has a list  of colours $\calc_v$ that may be assigned to $v$ of size 
\[q:=200r\frac{\D\ln\ln\D}{\ln\D}.\]
A partial list colouring $\s$ is an assignment to a subset of the vertices, where the colours are drawn from their lists.  Given any partial colouring, $L_v$ is defined to the the set of colours in $\calc_v$ not appearing on any neighbours of $v$ along with $\blank$.}

Because we are not trying for a good constant, we can afford to be a bit looser in our definition of $L$ and our second flaw will be simpler than that in Section~\ref{stf}. We define

\[L=\D^{9/10}.\]

Given a partial colouring $\s$, we define the following two {\em flaws} for any vertex $v$:

\begin{eqnarray*}
B_v&\equiv & |L_v|<L\\
Z_v &\equiv& \mbox{ at least $L$ neighbours of $v$ are coloured } \blank.
\end{eqnarray*}

{
\begin{observation}\label{odist2} $B_v$ and $Z_v$ are determined by the colours of the vertices in $N(v)$.
\end{observation}
}

It is trivial to see that any flaw-free partial colouring can be completed greedily to a full colouring of $G$,
as the list of available colours for each vertex is greater than the number of uncoloured neighbours.

Again, we say $v$ is the {\em vertex} of flaw $f=B_v$ or $Z_v$, and we denote $v(f):=v$. We use the same ordering on the flaws:
Every $B_v$ comes before every $Z_u$, and the $B_v$'s and $Z_u$'s are each ordered according to the labels of $v,u$.

We find a flaw-free partial colouring using essentially the same algorithm we used for triangle-free graphs, but we must be more careful about recolouring a neighbourhood.  It will be useful to represent a partial colouring of a neighbourhood as a collection of disjoint independent sets.

{We let $\calc=\cup_{v\in G}\calc_v$ denote the set of all colours that may appear in the graph, and define:}

\begin{definition}\label{dpc} {Given a vertex $v$ and a fixed partial colouring of $V(G)\bk N_v$,} a {\em partial colour assignment} to $N_v$ is a collection of disjoint independent sets $(\theta_1,...,\theta_{|\calc|})$, each a subset of $N_v$, such that for any $u\in \theta_i$ we have:
$i\in\calc_u$ and $i$ does not appear on any neighbour of $u$ outside of $N_v$.
\end{definition}

It is possible that $\theta_i=\emptyset$, and we do not require that $\cup_{i=1}^{|\calc|}\theta_i=N_v$.  Any $u\in N_v$ that is not in any of the $\theta_i$ is considered to be coloured $\blank$.

To recolour $N_v$, we take a uniformly random partial colour assignment to $N_v$ and then assign the colour $i$ to every vertex in each $\theta_i$.  More specifically, given a colouring $\s$ and a vertex $v$, we let $\Omega$ denote the set of all partial colour assignments to $N_v$ and we choose a uniform member of $\Omega$.

Note that if $N_v$ contains no edges, then this recolouring is equivalent to giving each $u\in N_v$ a uniform colour from $N_u$, as we did in \fix.

We use the same flaw ordering as in Section~\ref{stf}; i.e.\ every $B_v$ comes before every $Z_u$, and the $B_v$'s and $Z_u$'s are each ordered according to the labels of $v,u$.  

{The following procedure differs from FIX only in the distances: Observation~\ref{odist2} allows us to recurse on flaws $Z_w$
within distance two rather than three.  And we increase the distance for flaws $B_w$ from two to three so that we get Observation~\ref{ozb} below, which will be very useful in our analysis.}

\begin{tabbing}
{\bf FIX2}($f,\s$)\\
Set $v=v(f)$.\\
Choose a uniformly random partial colour assignment to $N_v$ and then recolour $N_v$ accordingly.\\
While \= there are any flaws $B_w$ with $\dist(w,v)\leq 3$ or $Z_w$ with $\dist(w,v)\leq 2$:\\
\>Let $g$ be the least such flaw and call {\bf FIX}($g,\s'$) where $\s'$ is the current colouring.\\
Return the current colouring.
\end{tabbing}

{
\begin{observation}\label{ozb}
Whenever we call \fix$(Z_u,\s)$ we have that $B_w$ does not hold for any $w\in N_u$. 
\end{observation}

This observation follows from our flaw ordering, and the fact that we call \fix on flaws $B_w$ with $w$ up to distance three from $v$ rather than two.  
}

The analog of Observation~\ref{ofix} holds again here, and so to prove Theorem~\ref{mt2} it suffices to prove that \fixx terminates with positive probability.

We will assume throughout the remainder of this section that $\D\geq 2^{200r}$ as otherwise the bound of Theorem~\ref{mt2} is trivial.

\subsection{More probability bounds}  

We begin with some key lemmas from Shearer's paper on the independence number of a $K_r$-free graph\cite{js}.
We rephrase the short  proofs here for completeness and to extract a useful fact from them.

Given a graph $H$, we define:

\[ I(H) \mbox{ is the number of independent sets of $H$}.\]

\begin{lemma}\label{libig} {For any $r\geq 2$}, if $H$ is $K_r$-free then $2^{|V(H)|}\geq I(H)\geq2^{|V(H)|^{\frac{1}{r-1}}-1}$.
\end{lemma}

\proofstart The upper bound is simply the number of subsets of $V(H)$. For the lower bound, we will prove that $H$ has an independent set of size at least $|V(H)|^{1/{r-1}}-1$; the bound follows by considering all subsets of that independent set.

We proceed by induction on $r$.  The trivial base case is $r=2$.  For $r\geq 3$: If any vertex $u\in H$ has degree at least $d=|V(H)|^{\frac{r-2}{r-1}}$ then {since the neighbourhood of $u$ in $H$ is $K_{r-1}$-free, there is a sufficiently large independent set in that neighbourhood by induction.} Otherwise, the {maximum degree in $H$ is less than $d$ and so the} straightforward greedy algorithm finds an independent set of size at least $|V(H)|/(d+1)>|V(H)|^{1/{r-1}}-1$.
\proofend

\begin{lemma}\label{lmu} If $H\neq\emptyset$  is $K_r$-free, $r\geq4$, then half of the independent sets in $H$ have size at least $\inv{2r} \log_2 I(H)/\log_2\log_2 I(H)$.
\end{lemma}

\proofstart  It suffices to show that at most $\hf I(H)$ subsets of $V(H)$ have size at most $\ell=\lfloor \inv{2r} \log_2 I(H)/\log_2\log_2 I(H)\rfloor$; i.e:
\begin{equation}\label{e33}
\sum_{i=0}^{\ell}{|V(H)|\choose i}\leq\hf I(H).
\end{equation}
We can assume $\log_2 I(H)\geq2$ as otherwise $\ell=0$  and so the lemma is trivial (since $H\neq\emptyset$). We can also assume $r\leq \log_2 I(H)/2\log_2\log_2 I(H)$ else $\ell=0$. We set $x=\log_2 I(H)\geq 2$.  
Rearranging the second inequality of Lemma~\ref{libig} gives $|V(H)|\leq(1+\log_2 I(H))^{r-1}$ and so we substitute $h=(1+\log_2 I(H))^{r-1}\geq 27$ for $|V(H)|$ in~(\ref{e33}).  So $h=(1+x)^{r-1}<\inv{4} x^{2r}$ for $x\geq 2$. {Also, a simple induction on $\ell$ confirms that $\sum_{i=0}^{\ell}{h\choose i}\leq \sum_{i=0}^{\ell-1}{h\choose i} +\frac{h^{\ell}}{\ell!}\leq 2h^{\ell}$ for $\ell\geq0,h\geq 2$.  }  
So the LHS of~(\ref{e33}) is at most 

\[ 2h^{\ell} <\hf x^{2r\ell}\leq\hf 2^{\log_2 x\times\frac{x}{\log_2 x}}
=\hf 2^x=\hf I(H).\]
 This proves~(\ref{e33}).
 \proofend

\noindent {\bf Remarks} \\
\noindent(1) Lemma~\ref{libig} is the only place where we use the fact that our graph is $K_r$-free.  Our proof shows that the bound of Theorem~\ref{mt2} holds whenever every subgraph $H\subseteq G$ satisfies the implication of either Lemma~\ref{libig} or Lemma~\ref{lmu}.  In fact, it is enough for this to hold for every $v$ and $H\subseteq N(v)$.

\noindent(2) Note that the argument in Lemma~\ref{lmu} can in fact show that the average size of the independent sets of $H$ is at least $\inv{2r} \log_2 I(H)/\log_2\log_2 I(H)$, which is Lemma 1 of \cite{js}.  

\noindent(3)  Alon~\cite{na} proves that if $G$ is locally $r$-colourable, meaning that every neighbourhood can be $r$-coloured, then for any $v$ and $H\subseteq N_v$, the median size of the independent sets of $H$ is at least $\inv{10\log_2(r+1)} \log_2 I(H)$.  Plugging this bound into the rest of our proof yields that  $\chi_{\ell}\leq O(\ln r \frac{\D}{\ln\D})$ for such graphs, as shown in~\cite{aj2}.

We use these to bound the probabilities of our flaws.

{\bf Setup for Lemma~\ref{lbv2}:}  Each vertex $u\in N_v$ has a list {$L^*_u$ containing $\blank$ and perhaps other colours; specifically, $L^*_u$ is the set of colours of $\calc_u$ not appearing on any neighbour of $u$ outside of $N_v$ along with $\blank$.} We give the vertices of $N_v$ a random partial colour assignment consistent with these lists. This assignment determines $L_v$ - the set of colours in $\calc_v$ that do not appear in the partial colour assignment.  

\begin{lemma}\label{lbv2}  
\begin{enumerate}
\item[(a)] $\pr(|L_v|<L)<\D^{-4}$.
\item[(b)] The probability that at least $L$ neighbours of $v$ are coloured $\blank$ and $|L_u|>L$ for all $u\in N_v$ is at most $\D^{-4}$.
\end{enumerate}
\end{lemma}

\proofstart  We begin with a method for sampling a  partial colour assignment.  

Define $\Omega$ to be the set of all partial colour assignments to $N_v$, and let $W=(W_1,...,W_{|\calc|})$ be a uniform member of $\Omega$. Define $Q_1$ to be the vertex set consisting of $W_1$ and all blank vertices which can be given the colour 1; {i.e.\ all blank  $u\in N_v$ with $1\in L^*_u$.}   Select a uniformly random independent set $W_1'$ of $Q_1$ and form $W'$ by replacing $W_1$ with $W_1'$.

{\em Claim 1:  $W'$ is a uniform member of $\Omega$.} 

{\em Proof of Claim 1:} For any $|\calc| -1$ disjoint independent sets $S_2,...,S_{|\calc|}\subseteq N_v$ we define
$\Omega_{S_2,...,S_{|\calc|}}\subseteq\Omega$ to be the set of partial colour assignments $(\theta_1,...,\theta_{|\calc|})$ with $\theta_2=S_2,...,\theta_{|\calc|}=S_{|\calc|}$; so this yields a partition of $\Omega$.  Note that $W'$ is a uniform member of $\Omega_{W_2,...,W_{|\calc|}}$.  Furthermore, because $W$ is a uniform member of $\Omega$, the part $\Omega_{W_2,...,W_{|\calc|}}$ is selected with the correct distribution, i.e with probability $|\Omega_{W_2,...,W_{|\calc|}}|/|\Omega|$.  So $W'$ is a uniform member of $\Omega$.  \proofend

Repeating this argument, we can resample $W_2,...,W_{|\calc|}$ in the same manner. Specifically:

\begin{tabbing}
Let  $W=(W_1,...,W_{|\calc|})$ be a uniform member of $\Omega$.\\
For \=$i=1$ to $|\calc|$\\
\>Define $Q_i$ to be the subgraph induced by $W_i$ and \= all vertices {that are blank at this step}\\ 
\>\>and can be given the colour $i$.\\
\>Let $W_i'$ be a uniform independent set of $Q_i$\\
\>Modify $W$ by replacing $W_i$ with $W_i'$.
\end{tabbing}

{To be clear: the blank vertices in the definition of $Q_i$ are blank in the current partial colour assignment} $W=(W'_1,...,W_{i-1}',W_i,...,W_{|\calc|})$.  By repeating the argument from Claim 1, we see that the partial colour assignment produced by this procedure is a uniform member of $\Omega$.  

{\em Part (a):}  Let $A_1$ be the set of colours $i\in \calc_v$ such that $I(Q_i)\leq \D^{1/20}$, and set $A_2 := \calc_v\bk A_1$. Since the subgraph induced by $N_v$ is $K_{r-1}$-free, Lemma~\ref{lmu} implies that  for each $i\in A_2$  the median independent set of $Q_i$ has size at least $\inv{2(r-1)}\log_2 I(Q_i)/\log_2\log_2 I(Q_i)>\inv{40r}\log_2\D/\log_2\log_2\D$. (When applying Lemma~\ref{lmu} note that if $Q_i=\emptyset$ then $i\in A_1$.)

{At iteration $i$: If colour $i\in A_1$ then the probability that we choose $W_i'=\emptyset$} is $\inv{I(Q_i)}\geq\D^{-1/20}$.  Note that if $W_i'=\emptyset$ then $i$ will be in $L_v$.  {If $i\in A_2$, then with probability at least $\hf$ we choose a $W_i'$ with} $|W_i'|\geq \inv{40r}\log_2\D/\log_2\log_2\D$.  Since the total size of the sets $W_i'$ is at most $\D$, this cannot happen for more than $40r\frac{\D\log_2\log_2\D}{\log_2\D}$ colours.

We consider two random binary strings, each of length $|\calc_v|$.  In the first, each bit is 1 with probability $\D^{-1/20}$, and 0 otherwise.  In the second, the bits are uniform.  By coupling the choice of $W_i'$ with these bits, we ensure that: (a) for each $i\in A_1$, if the corresponding bit in the first stream is 1 then $W_i'=\emptyset$; (b) for each $i\in A_2$, if the  corresponding bit in the second stream is 1 then $|W_i'|\geq  \inv{40r}\log_2\D/\log_2\log_2\D$.  For example, in iteration $i$ if we have $I(Q_i)<\D^{1/20}$ and so $i\in A_1$ then we look at the next bit of the first string.  If that bit is 1 then we set $W_i'=\emptyset$; otherwise we set $W_i'=\emptyset$ with probability $\frac{1}{I(Q_i)}-\D^{-1/20}$. Similarly when $i\in A_2$.

{
Set $\ell=\hf |\calc_v| =100r\D\log_2\log_2\D/\log_2\D$, and so we must have either $A_1\geq\ell$ or $|A_2|\geq\ell$.

{\em Claim 2: If  the outcomes of this procedure yield $|L_v|<L$ then at least one of these two events must hold:
\begin{itemize}
\item  $E_1 =$ at most $L$ of the first $\ell$ bits of the first string are 1 
\item $E_2=$ at most $40r\frac{\D\log_2\log_2\D}{\log_2\D}$ of the first $\ell$ bits of the second stream are 1 
\end{itemize}
}

{\em Proof:} If $W'_i=\emptyset$ then $i\in L_v$.  So $\ov{E_1}$ and the event $|A_1|\geq\ell$ imply that at least $L$ colours in $A_1$ are in $L_v$. $\ov{E_2}$ and the event $|A_2|\geq\ell$ imply that for more than $40r\frac{\D\log_2\log_2\D}{\log_2\D}$ colours $i\in A_2$ we have $|W_i'|\geq \inv{40r}\log_2\D/\log_2\log_2\D$, which contradicts the fact that the sets $W_i'$ are disjoint and have total size at most $|N_v|\leq\D$.  Since we must have either 
$|A_1|\geq\ell$ or $|A_2|\geq\ell$  then if $|L_v|<L$ we must have $E_1\vee E_2$. 
\proofend

Claim 2 implies $\pr(|L_v| < L)\leq\pr(E_1)+\pr(E_2)$.}
Note that the expected number of 1's in the first $\ell$ bits of the first string is $\ell\times \D^{-1/20}\gg L=\D^{9/10}$ and the expected number of 1's in the first $\ell$ bits of the second string is $\hf\ell = 50r\D\log\log\D/\log\D$. So  the Chernoff Bounds (or Lemma~\ref{lncv}) imply that each of $E_1,E_2$ occur with probabilty less than $\hf\D^{-4}$ for $r\geq 4$ and $\D\geq 2^{500 r}$.  This proves part (a).

{\em Part (b):}  Consider any $L$ neighbours $u_1,...,u_L\in N_v$.  We will prove the probability that each $u_i$ is coloured blank and satisfies $|L_{u_i}|>L$ is at most $1/L!$.  This proves part (b) as ${\D\choose L}/L!<\D^{-4}$ for $\D\geq 100$.

{Fix a colouring of $V(G)\bk N_v$ and}
let $\Omega_B\subset\Omega$ be the set of partial colour assignments in which every $u_i$ is coloured $\blank$ and satisfies $|L_{u_i}|>L$. (Note: a partial colour assignment in $\Omega_B$ may also have additional blank vertices.) Take any $W\in \Omega_B$ and extend it to a partial colour assignment $W_2$ in which each of $u_1,...,u_t$ are not blank as follows:

\begin{tabbing}
begin with the colouring $W$\\
for \= $i=1$ to $L$\\
\>give $u_i$ a colour from $L^*_u$ which does not appear on any of its neighbours in $N_v$.
\end{tabbing}

This yields a colouring $W'$ of $N_v$ which can be viewed as the partial colour assignment $(\theta_1,...,\theta_{|\calc|})$ where $\theta_j$ is the set of vertices with colour $j$ in $W'$.

By definition of $\Omega_B$, each $u_i$ has at least $L$ available colours in $W$.  By the time we reach iteration $i$, at most $i-1$ of those colours have been assigned to a neighbour of $u_i$ in $\{u_1,...,u_{i-1}\}$.  So there are always at least $L-i+1$ choices for a colour to assign to $u_i$ and so the number of choices for $W'$ is at least $L!$.
Each partial colour assignment $W'$ can arise from at most one $W\in\Omega_B$, namely the $W$ obtained from $W'$ by colouring $u_1,...,u_L$ all $\blank$.  So $|\Omega_B|\leq |\Omega|/L!$, which is what we need to establish part (b).
\proofend

\subsection{FIX2 terminates}
Now the same argument from Section~\ref{sat} implies that  FIX2 terminates with positive probability, and thus proves Theorem~\ref{mt2}.  

Each time we call \fixx$(v,\s)$  we let $\call=\{L^*_u: u\in N_v\}$ be the lists of available colours on the neighbours of $v$ in the colouring obtained from $\s$ by uncolouring $N_v$; i.e.\ $L^*_u$ is the set of colours in $\calc_u$ that do not appear on any neighbours of $u$ outside of $N_v$, along with $\blank$.  We let $\Omega(\call)$ be the set of partial colour assignments to $N_v$ consistent with $\call$.  We let $\calb(\call)\subset\Omega(\call)$ be the set of partial colour assignments that have the flaw $B_v$.  We let $\calz(\call)\subset\Omega(\call)$ be the set of partial colour assignments which have the flaw $Z_v$.

We define $H_t,R_t$ analogously to Section~\ref{sat}.  At each step: If we are addressing the flaw $B_v$ then Lemma~\ref{lbv2}(a) implies that the number of choices for the colouring of $N_v$ before the recolour line is at most $|\calb(\call)|\leq \D^{-4}|\Omega(\call)|$.  If we are addressing the flaw $Z_v$ then by Observation~\ref{ozb}, each $u\in N_v$ has at least $L$ available colours in $\s$ and so must have $|L_u|\geq L$ {before} uncolouring $N_v$; thus $|L^*_u|\geq|L_u|\geq L$.   So Lemma~\ref{lbv2}(b) implies that the number of choices for the colouring of $N_v$ before the recolour line is at most $|\calz(\call)|\leq \D^{-4}|\Omega(\call)|$.  This yields that the size of what is written to $H_t$ is $3\log_2\D+\log_2|\Omega(\call)|-4\log_2\D+O(1)$ whereas the number of random bits used is $\log_2|\Omega(\call)|$. This is enough for the analysis from Section~\ref{sat}, in particular the proof of Lemma~\ref{cshort} to carry through.  

\noindent {\bf Remark} This time it is not clear how to obtain a polytime algorithm; the challenge is  to select a uniform  partial colour assignment efficiently.  Johansson's proof yields a polytime algorithm (see~\cite{gg}).

{
\section{Lopsided Local Lemma}  Bernshteyn notes that the proofs of Theorems~\ref{mt} and~\ref{mt2} could have been carried out using the Lopsided Local Lemma rather than an entropy compression argument.  One considers taking a uniformly random partial colouring of the entire graph.  The bad events are: $B_v$ and $Z_v\wedge \ov{B_v}$. By conditioning on the colours of all vertices at distance at least two or three from $v$, Lemmas~\ref{lbv} and~\ref{lbv2} imply that the probability of the bad events is sufficiently small, even when conditioning on the outcomes of distant events.  See~\cite{ab} for more details and for an extension of these results to DP-colouring.
}

 \section*{Acknowledgement}  My thanks to Dimitris Achlioptas and Fotis Iliopoulos for some very helpful discussions. I am also grateful to two anonymous referees for many helpful comments and to Zdenek Dvorak for pointing out a problem with Lemma~\ref{lncv} in an earlier version. This research is supported by an NSERC Discovery grant.

\end{document}